\documentclass{amsart}
\pdfoutput=1 
\usepackage{amssymb, latexsym, amsmath, verbatim, amsthm, amscd}
\usepackage{mathtools}
\usepackage{caption}
\usepackage{pinlabel}
\usepackage{graphicx}
\usepackage{hyperref}

\title{Morse subsets of CAT(0) spaces are strongly contracting}
\author{Christopher H. Cashen} 
\address{
Faculty of Mathematics,
University of Vienna,
1090 Vienna \\Austria}
\email{\href{mailto:christopher.cashen@univie.ac.at}{christopher.cashen@univie.ac.at}}
\hypersetup{
    pdftitle={Morse subsets of CAT(0) spaces are strongly contracting},    % title
    pdfauthor={Christopher H. Cashen},     % author
    pdfkeywords={Morse set, contracting set, recurrent set, CAT(0)
      space, strongly contracting}, % list of keywords
    colorlinks=true,       % false: boxed links; true: colored links
    linkcolor=black,          % color of internal links
    citecolor=black,        % color of links to bibliography
    filecolor=black,      % color of file links
    urlcolor=black           % color of external links
}

\theoremstyle{plain}
\newtheorem*{theorem}{Theorem}
\newtheorem*{proposition}{Proposition}
\newtheorem*{corollary}{Corollary}

\mathtoolsset{centercolon} % makes the : and = in := be same height

\newcommand{\from}{\colon\thinspace} % colon for maps

\newcommand{\Z}{\mathcal{Z}}
\newcommand{\X}{\mathcal{X}}
\DeclareMathOperator{\diam}{diam}
\DeclareMathOperator{\len}{len}

\begin{document}
\begin{abstract}
  We prove that Morse subsets of CAT(0) spaces are strongly contracting.
  This generalizes and simplifies a result of Sultan, who
  proved it for Morse quasi-geodesics.
  Our proof goes through the recurrence characterization of Morse
  subsets. 
\end{abstract}
\maketitle
%---------------------------------------
% End of frontmatter
%---------------------------------------
%---------------------------------------
% Start of main body of article
%---------------------------------------

%!TEX root = main.tex

In this note we give a short proof of the following technical result:
\begin{proposition}\label{prop}
  If $\Z$ is a closed, $\rho$--recurrent subset of a CAT(0) space then
  $\Z$ is $12\rho(21)$--strongly contracting. 
\end{proposition}

This is the final piece of the following theorem, which says that a
number of properties that are equivalent to quasi-convexity in
hyperbolic spaces are also equivalent to one another in CAT(0) spaces:

\begin{theorem}
  Let $\X$ be a geodesic metric space. Let $\Z$ be a closed, unbounded subset of
  $\X$.
  For $x\in\X$, let $\pi_\Z(x):=\{z\in\Z\mid d(x,z)=d(x,\Z)\}$.
  The following are equivalent:
  \begin{description}
  \item [$\Z$ is Morse] There is a function $\mu\from [1,\infty)\times
    [0,\infty)\to [0,\infty)$ defined by
    $\mu(L,A):=\sup_\gamma\sup_{w\in\gamma} d(w,\Z)$, where the first
    supremum is taken over $(L,A)$--quasi-geodesic segments $\gamma$
    with both endpoints on $\Z$.

  \item [$\Z$ is contracting] There is a function $\sigma\from
    [0,\infty)\to[0,\infty)$ defined by
    $\sigma(r):=\sup_{d(x,y)\leq d(x,\Z)\leq r}\diam
    \pi_\Z(x)\cup\pi_Z(y)$ that satisfies
    $\lim_{r\to\infty}\sigma(r)/r=0$.
   \item [$\Z$ is recurrent] There is a function $\rho\from
      [1,\infty)\to[0,\infty)$ defined by 
      $\rho(q):=\sup_{\Delta(\gamma)\leq q}\inf_{w\in\gamma}
    d(w,\Z')$,
    where the first supremum is taken over rectifiable segments
    $\gamma$ with endpoints $z,\,z'\in\Z$ such that
    $\Delta(\gamma):=\frac{\len(\gamma)}{d(z,z')}\leq q$ and $\Z'$ is $\Z$
    with the open balls of radius $d(z,z')/3$ about $z$ and $z'$ removed. 
  \end{description}

   If $\X$ is hyperbolic or CAT(0) then these conditions are equivalent to:
     \begin{description}

   \item [$\Z$ is strongly contracting] $\Z$ is contracting and the
    contraction gauge $\sigma$ is a bounded function.
   \end{description}
 \end{theorem}
 
We refer the reader to \cite{BriHae99} for background on  hyperbolic and CAT(0) spaces.
  
 \begin{corollary}
   Morse subsets of CAT(0) spaces are strongly contracting.
 \end{corollary}

The corollary confirms a conjecture of Russell, Spriano, and Tran
\cite{RusSprTra18} and generalizes results of Sultan \cite{Sul14},
who proved that Morse quasi-geodesics in CAT(0) spaces are strongly
contracting, and Genevois \cite{Gen17}, who proved that Morse subsets of CAT(0) cube
complexes are strongly contracting.

The condition that $\Z$ is closed is inessential. 
It guarantees that
the empty set is not in the image of $\pi_\Z$. 
This hypothesis can be
avoided by defining $\pi_\Z(x):=\{z\in\Z\mid d(x,z)\leq d(x,\Z)+1\}$.
Extra bookkeeping is then required to compute an explicit contraction
bound in the proof of the proposition.

The four properties are trivially satisfied for bounded sets,
with the possible exception that recurrence can fail if some $\Z'$ is empty. 
For example, a two point set is not recurrent, but its contraction
gauge is bounded by its diameter.

\begin{proof}[Proof of the theorem]
  The contraction condition was introduced in \cite{ArzCasGrub}, where
  it was shown to be equivalent to the Morse condition. The recurrence
  condition was used to characterize Morse quasi-geodesics in
  \cite{DruMozSap10}, and this characterization can be extended to
  arbitrary subsets, as in \cite[Theorem~2.2]{CasMac17}.
  Strong contraction obviously implies contraction.
  It is easy to see that all of these properties are equivalent to
  quasi-convexity in hyperbolic spaces.
  The proposition supplies the remaining implication.
\end{proof}

  There is extensive literature making use of the Morse property and
  equivalent characterizations in various settings, but a complete
  exposition would be longer than this paper, so we will not attempt it.
  Sultan's result uses a characterization of the images of Morse
  quasi-geodesics in asymptotic cones due to Dru\c{t}u, Mozes, and Sapir \cite{DruMozSap10}.
  Loosely speaking, this characterization depends on there being a
  sensible notion of one point being \emph{between} two others, which we have
  for quasi-geodesics but not, at least in an obvious way, for
  arbitrary subsets.
  We avoid the use of asymptotic cones and instead use recurrence
  (which also comes from \cite{DruMozSap10}).
  We construct curves in essentially the same way as Sultan, but our argument,
  in addition to applying to general subsets, is simpler and
  gives an explicit strong contraction bound.

\begin{proof}[Proof of the proposition]
  Define $D:=\rho(21)$.
  Supposing the contraction gauge $\sigma$ of $\Z$ is not bounded by
  $12D$, we derive a contradiction.
 Failure of the contraction bound means there exist points $x,\,y\in \X$ such that
  $d(x,y)\leq d(x,\Z)$ and such that $\diam
  \pi_\Z(x)\cup\pi_\Z(y)>12D$.
  We may assume $d(x,\Z)\geq d(y,\Z)$, because otherwise $d(x,y)\leq
  d(y,\Z)$ and we can swap the roles of $x$ and $y$.
  Choose $x'\in\pi_\Z(x)$ and $y'\in\pi_\Z(y)$ such that $P:=d(x',y')>12D$.
  Let $\Z'$ denote the set $\Z$ with the open balls of radius $P/3$
  about $x'$ and $y'$ removed.

  For points $a,\,b\in \X$, let $[a,b]\from [0,1]\to\X$ denote
  the geodesic segment from $a$ to $b$, parameterized
  proportional to arc length.
  Concatenation is denoted `$+$'.  
  \begin{equation}
    \label{ob}
       \text{If $d(w,\Z')\leq D$ for some $w\in \X$ then $w\notin [x',x]+[x,y]+[y,y']$.}\tag{$*$}
  \end{equation}
  To see this, first suppose $w\in [x',x]$.
   Then $x'\in\pi_\Z(w)$, so $P/3\leq d(x',\Z')\leq
  d(x',w)+d(w,\Z')=d(w,\Z)+d(w,\Z')\leq 2d(w,\Z')\leq 2D$, which is a
  contradiction, since $P>12D$.
  Similarly, $w\notin [y',y]$.
  If $w\in [x,y]$ then:
  \[d(x,w)+d(w,y)=d(x,y)\leq d(x,\Z)\leq d(x,w)+D\]
  Thus, $d(w,y)\leq D$, which implies:
  \[P/3\leq d(y',\Z')\leq d(y',y)+d(y,\Z')\leq 2d(y,\Z')\leq
    2(d(y,w)+d(w,\Z'))\leq 4D\]
  Again, this contradicts the hypothesis that $P>12D$, so \eqref{ob}
  is verified.

  \par{\bf Case 1, $d(x,x')\leq 6P$: }
  Define $\gamma:=[x',x]+[x,y]+[y,y']$.
  Then $\len(\gamma)\leq 18P<21P$, so recurrence says there
  is a point $w\in \gamma$ such that $d(w,\Z')\leq D$.
  By \eqref{ob}, this is impossible.

\begin{figure}[h]
    \begin{minipage}[b]{0.3\linewidth}
      \labellist
      \small
      \pinlabel $x$ [r] at 2 135
      \pinlabel $x'$ [r] at 3 2
      \pinlabel $y$ [l] at 80 90
      \pinlabel $y'$ [l] at 80 2
      \endlabellist
\centering
\includegraphics[height=.9in]{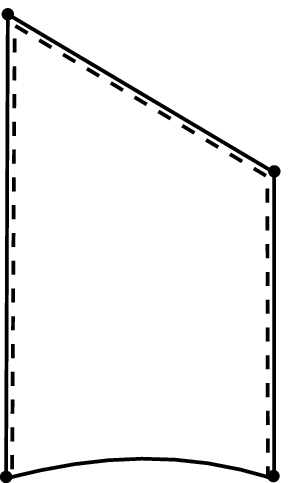}
\caption*{Case 1}
\end{minipage}
\hspace{0.3cm}
\begin{minipage}[b]{0.3\linewidth}
      \labellist
      \small
      \pinlabel $x$ [r] at 2 135
      \pinlabel $x'$ [r] at 3 2
      \pinlabel $y$ [l] at 80 90
      \pinlabel $y'$ [l] at 80 2
      \pinlabel $a$ [r] at 0 70
      \pinlabel $b$ [bl] at 40 115
      \endlabellist
\centering
\includegraphics[height=.9in]{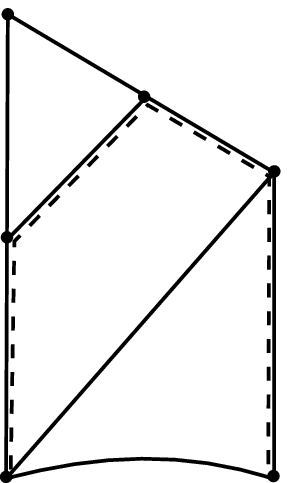}
\caption*{Case 2}
\end{minipage}
\hspace{0.3cm}
\begin{minipage}[b]{0.3\linewidth}
        \labellist
      \small
      \pinlabel $x$ [r] at 2 135
      \pinlabel $x'$ [r] at 3 2
      \pinlabel $y$ [l] at 80 90
      \pinlabel $y'$ [l] at 80 2
      \pinlabel $a$ [r] at 2 52
      \pinlabel $c$ [bl] at 35 55
      \pinlabel $e$ [tr] at 46 44
      \pinlabel $b$ [l] at 80 58
      \endlabellist
\centering
\includegraphics[height=.9in]{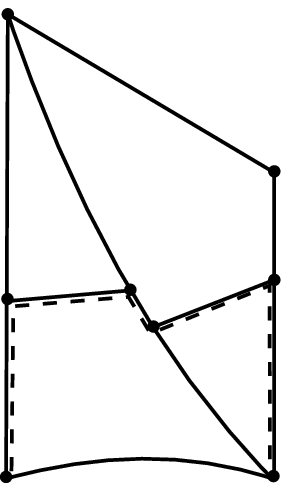}
\caption*{Case 3}
\end{minipage}
  \end{figure}
  
  \par{\bf Case 2, $d(x,x')>6P$ and $d(y,y')\leq 4P$: }
  Let $a:=[x',x](\frac{6P}{d(x,x')})$ and
  $b:=[y,x](\frac{6P}{d(x,x')})$, so that $d(a,x')=\frac{6P}{d(x,x')}\cdot d(x,x') =6P$ and
  $d(b,y)=\frac{6P}{d(x,x')}\cdot d(x,y)\leq 6P$.
  Since $d(x',y)\leq 5P$, the CAT(0) condition implies $d(a,b)<
  5P$.
  Define $\gamma:=[x',a]+[a,b]+[b,y]+[y,y']$.
  Since $\len(\gamma)\leq 6P+5P+6P+4P=21P$,  recurrence
  says there is a point $w\in\gamma$ with $d(w,\Z')\leq D$.
  By \eqref{ob}, $w\in [a,b]$, but this is impossible because
  $d([a,b],\Z)\geq d(a,\Z)-d(a,b)> 6P-5P=P>D$.

  \par{\bf Case 3, $d(x,x')>6P$ and $d(y,y')> 4P$: }
  Let $a:= [x',x](\frac{4P}{d(x,x')})$ and
let  $c:=[y',x](\frac{4P}{d(x,x')})$.
Then $d(x',a)=4P$ and:
\[4P\leq d(y',c)=\frac{4P}{d(x,x')}\cdot d(y',x)\leq
\frac{4P}{d(x,x')}\cdot (d(x,x')+P)\leq \frac{14}{3}P\]

  Let $b$ be the point of $[y',y]$ at
  distance $4P$ from $y'$, and let $e$ be the point of $[y',x]$ at
  distance $4P$ from $y'$, so
$d(c,e)\leq \frac{2}{3}P$.
The CAT(0) condition implies that $d(a,c)< P$ and, since $d(x,y)\leq
d(x,y')$, that $d(e,b)\leq 4\sqrt{2}P$.

Define $\gamma:=[x',a]+[a,c]+[c,e]+[e,b]+[b,y']$. 
Then $\len(\gamma)< 4P+P+\frac{2}{3}P+4\sqrt{2}P+4P<21P$, so
recurrence demands a point $w\in\gamma$ with $d(w,\Z')\leq D$.
By \eqref{ob}, $w\notin [x',a],\, [b,y']$.
We cannot have $w\in [a,c]+[c,e]$ because $d([a,c]+[c,e],\Z)\geq
d(a,\Z)-(d(a,c)+d(c,e))\geq 4P-P-\frac{2}{3}P>D$.
Thus, $w\in [e,b]$, so $d(e,b)=d(e,w)+d(w,b)$.
However, $d(w,b)\geq d(b,\Z)-d(w,\Z)\geq 4P-D>\frac{47}{12}P$.
By the same reasoning, $\frac{47}{12}P<d(a,w)$, but $d(a,w)\leq
P+\frac{2}{3}P+d(e,w)$, so $d(e,w)>\frac{27}{12}P$.
This gives us a contradiction:
\[6P<\frac{74}{12}P<d(e,w)+d(w,b)=d(e,b)\leq 4\sqrt{2}P<6P\qedhere\]
\end{proof}

%---------------------------------------
% End of main body of article
%---------------------------------------
%---------------------------------------
% Start of backmatter
%---------------------------------------

%---------------------------------------
% End of backmatter
%---------------------------------------

\end{document}